\documentclass[12pt]{article}

\usepackage{amsmath, epsfig, cite}
\usepackage{amssymb}
\usepackage{amsfonts}
\usepackage{latexsym}
\usepackage{float}
\usepackage{color}

\newtheorem{thm}{Theorem}[section]
\newtheorem{prop}[thm]{Proposition}

\newtheorem{lem}[thm]{Lemma}

\numberwithin{equation}{section}

\newcommand{\qed}{{\hfill$\square$}\medskip}
\begin{document}

\begin{center}
{\Large\bf On the divisibility of $q$-trinomial coefficients}
\end{center}

\vskip 2mm \centerline{Ji-Cai Liu}
\begin{center}
{\footnotesize Department of Mathematics, Wenzhou University, Wenzhou 325035, PR China\\
{\tt jcliu2016@gmail.com } \\[10pt]
}
\end{center}


\vskip 0.7cm \noindent{\bf Abstract.}
We establish a congruence on sums of central $q$-binomial coefficients. From this $q$-congruence, we derive the divisibility of the $q$-trinomial coefficients introduced by Andrews and Baxter.


\vskip 2mm
\noindent{\it MR Subject Classifications}: 05A10, 11A07, 11B65, 13A05

\section{Introduction}
We consider the coefficients of the expanded form of this expression:
\begin{align}
(1+x+x^2)^n=\sum_{j=-n}^{n}\left( n \choose j \right)x^{j+n}.\label{new-1}
\end{align}
The coefficients in \eqref{new-1} are called trinomial coefficients.
Two simple formulas for the trinomial coefficients (see \cite[page 43]{sills-b-2018}) are
\begin{align*}
\left( n \choose j \right)=\sum_{k=0}^n{n\choose k}{n-k\choose k+j},
\end{align*}
and
\begin{align*}
\left( n \choose j \right)=\sum_{k=0}^n(-1)^k{n\choose k}{2n-2k\choose n-j-k}.
\end{align*}

In 1987, Andrews and Baxter \cite{ab-1987} introduced six apparently distinct $q$-analogues of
the trinomial coefficients, such as
\begin{align}
\left({n\choose j}\right)_q=\sum_{k=0}^{n}
q^{k(k+j)}{n \brack k} { n-k \brack k+j},\label{new-2}
\end{align}
where the $q$-binomial coefficients are defined as
$$
{n\brack k}={n\brack k}_q
=\begin{cases}
\displaystyle\frac{(1-q^n)(1-q^{n-1})\cdots (1-q^{n-k+1})}{(1-q)(1-q^2)\cdots(1-q^k)}, &\text{if $0\leqslant k\leqslant n$},\\[10pt]
0,&\text{otherwise.}
\end{cases}
$$
These $q$-trinomial coefficients play an important role in the solution of a model in statistical mechanics (see \cite{ab-1987}). In 1990, Andrews \cite{andrews-jams-1990} investigated some interesting properties for the $q$-trinomial coefficients, which leads him to finite versions of dissections of the Rogers-Ramanujan identities into even and odd parts.

In 1999, Andrews \cite{andrews-dm-1999} showed that Babbage's congruence \cite{babbage-epj-1819}:
\begin{align}
{2p-1\choose p-1}\equiv 1\pmod{p^2}\label{new-3}
\end{align}
possesses the following nice $q$-analogue:
\begin{align}
{2p-1\brack p-1}\equiv q^{\frac{p(p-1)}{2}}\pmod{[p]_q^2},\label{new-4}
\end{align}
for any odd prime $p$. Notice that Wolstenholme \cite{wolstenholme-qjpam-1862} proved that \eqref{new-3} holds modulo $p^3$ for any prime $p\ge 5$, which is known as the famous Wolstenholme's theorem.

In order to understand \eqref{new-4}, we recall some necessary notation.
For polynomials $A_1(q), A_2(q),P(q)\in \mathbb{Z}[q]$, the $q$-congruence $$A_1(q)/A_2(q)\equiv 0\pmod{P(q)}$$ is understood as $A_1(q)$ is divisible by $P(q)$ and $A_2(q)$ is coprime with $P(q)$. In general, for rational functions $A(q),B(q)\in \mathbb{Z}(q)$,
\begin{align*}
A(q)\equiv B(q)\pmod{P(q)}\Longleftrightarrow
A(q)-B(q)\equiv 0\pmod{P(q)}.
\end{align*}
The $q$-integers are defined as $[n]_q=(1-q^n)/(1-q)$ for $n\ge 1$, and the $n$th
cyclotomic polynomial is given by
\begin{align*}
\Phi_n(q)=\prod_{\substack{1\le k \le n\\[3pt](n,k)=1}}
(q-\zeta^k),
\end{align*}
where $\zeta$ denotes an $n$th primitive root of unity.

Recently, Straub \cite[Theorem 2.2]{straub-pams-2019} established a $q$-analogue of Wolstenholme--Ljunggren congruence:
\begin{align*}
{ap\choose bp}\equiv {a\choose b}\pmod{p^3},
\end{align*}
for any prime $p\ge 5$, i.e.
\begin{align}
{an\brack bn}\equiv {a\brack b}_{q^{n^2}}-(a-b)b{a\choose b}\frac{n^2-1}{24}(q^n-1)^2\pmod{\Phi_n(q)^3}.
\label{new-5}
\end{align}
The simple case $a=2$ and $b=1$ in \eqref{new-5} reads
\begin{align}
{2n\brack n}\equiv 1-q^{2n^2}-\frac{n^2-1}{12}(q^n-1)^2\pmod{\Phi_n(q)^3},\label{a-3}
\end{align}
which is an extension of \eqref{new-4} (notice that ${2n-1\brack n}={2n\brack n}/(1+q^n)$).
We remark that Pan \cite[Lemma 3.1]{pan-mm-2013} obtained another interesting $q$-analogue of Wolstenholme--Ljunggren congruence and Zudilin \cite{zudilin-ac-2019} extended both Straub and Pan's $q$-congruences.

In the past few years, congruences for $q$-binomial coefficients as well as basic hypergeometric series attracted many experts' attention (see, for example, \cite{bachraoui-rj-2021,goro-ijnt-2019,guo-aam-2020,gs-rm-2020,gz-am-2019,lw-2020,liu-crm-2020,lp-a-2020,tauraso-aam-2012,WY2,wei-jcta-2021}).
Nowdays, there is also an extensive literature on congruences for trinomial coefficients (see, for instance, \cite{al-integers-2020,np-amc-2018,ws-a-2019,sunzw-scm-2014}).
It is worth mentioning that Gorodetsky \cite{goro-ijnt-2019} investigated congruences for the $q$-trinomial coefficients $\left({n\choose 0}\right)_q$.
However, the literature is still sparse on congruences for $q$-trinomial coefficients.

The first aim of the paper is to prove the following congruence for the $q$-trinomial coefficients $\left({n\choose 0}\right)_q$.
\begin{thm}\label{t-1}
For any positive integer $n$, the following congruence holds modulo $\Phi_n(q)^2$:
\begin{align*}
\left({n\choose 0}\right)_q\equiv {\begin{cases}
\displaystyle (-1)^m(1+q^m)q^{m(3m-1)/2}, &\text{if}~~n=3m,\\[10pt]
(-1)^mq^{m(3m+1)/2}, &\text{if}~~n=3m+1,\\[10pt]
(-1)^mq^{m(3m-1)/2},&\text{if}~~n=3m-1.
\end{cases}}
\end{align*}
\end{thm}

The elegant $q$-congruence \eqref{a-3} motivates us to establish the following congruence for the $q$-trinomial coefficients $({2n\choose n})_q$, which generalizes a conjectural $q$-congruence due to
Apagodu and the author \cite[Conjecture 1]{al-integers-2020}.
\begin{thm}\label{t-2}
For any positive integer $n$, the following congruence holds modulo $\Phi_n(q)^2$:
\begin{align*}
\left(2n\choose n\right)_q \equiv {\begin{cases}
\displaystyle 2(-1)^m(1+q^m)q^{m(3m-1)/2}-3m(1-q^{3m}), &\text{if}~~n=3m,\\[10pt]
2(-1)^mq^{m(3m+1)/2}-(3m+1)(1-q^{3m+1}), &\text{if}~~n=3m+1,\\[10pt]
2(-1)^mq^{m(3m-1)/2}-(3m-1)(1-q^{3m-1}),&\text{if}~~n=3m-1.
\end{cases}}
\end{align*}
\end{thm}

The rest of the paper is organized as follows. In Section 2, we shall mainly establish an auxiliary congruence on sums of central $q$-binomial coefficients, which is interesting by itself. The proofs of Theorems \ref{t-1} and \ref{t-2} will be presented in Sections 3 and 4, respectively.

\section{An auxiliary result}
In order to prove Theorems \ref{t-1} and \ref{t-2}, we first recall the following two lemmas.
\begin{lem} (See \cite[Lemma 2.3]{liu-cmj-2017}.)
For any non-negative integer $n$, we have
\begin{align}
&(1-q^n)\sum_{k=0}^{\lfloor n/2\rfloor}\frac{(-1)^kq^{k(k-1)/2}}{1-q^{n-k}}{n-k\brack k}\notag\\[10pt]
&={\begin{cases}
\displaystyle (-1)^m(1+q^m)q^{m(3m-1)/2}, &\text{if}~~n=3m,\\[10pt]
(-1)^mq^{m(3m+1)/2}, &\text{if}~~n=3m+1,\\[10pt]
(-1)^mq^{m(3m-1)/2},&\text{if}~~n=3m-1,
\end{cases}}\label{b-1}
\end{align}
where $\lfloor x \rfloor$ denotes the integral part of real $x$.
\end{lem}

\begin{lem} (See \cite[Lemma 3.3]{tauraso-aam-2012}.)
For $k=1,\cdots,n-1$, we have
\begin{align}
{2k-1\brack k}\equiv (-1)^kq^{k(3k-1)/2}{n-k\brack k} \pmod{\Phi_n(q)}.\label{b-2}
\end{align}
\end{lem}

We also require the following congruence on sums of the central $q$-binomial coefficients ${2k\brack k}$, which is interesting by itself.
\begin{prop}\label{p-1}
For any positive integer $n$, we have
\begin{align}
\sum_{k=1}^{\lfloor n/2 \rfloor}\frac{q^{-k(k-1)}}{[2k]_q} {2k\brack k}\equiv \frac{(1-q)(1-\mathcal{R}_n(q))}{1-q^n}\pmod{\Phi_n(q)},
\label{a-1}
\end{align}
where $\mathcal{R}_n(q)$ denotes the right-hand side of \eqref{b-1}.
\end{prop}

{\noindent\bf Remark.} Sun \cite[Theorem 1.1]{sunzw-scm-2011} proved that for any prime $p\ge 5$,
\begin{align}
\sum_{k=1}^{(p-1)/2}\frac{{2k\choose k}}{k}\equiv (-1)^{(p+1)/2}\frac{8p}{3}E_{p-3}\pmod{p^2},\label{new-6}
\end{align}
where $E_0,E_1,E_2,\cdots$ are Euler numbers. Notice that \eqref{a-1} is a $q$-analogue of \eqref{new-6} (the modulo $p$ version).

{\noindent\it Proof of \eqref{a-1}.}
By \eqref{b-2} and $q^n\equiv 1\pmod{\Phi_n(q)}$, we have
\begin{align}
\sum_{k=1}^{\lfloor n/2\rfloor}\frac{q^{-k(k-1)}}{[2k]_q} {2k\brack k}&=(1-q)\sum_{k=1}^{\lfloor n/2\rfloor}\frac{q^{-k(k-1)}}{1-q^{2k}} {2k\brack k}\notag\\[10pt]
&=(1-q)\sum_{k=1}^{\lfloor n/2\rfloor}\frac{q^{-k(k-1)}}{1-q^{k}} {2k-1\brack k}\notag\\[10pt]
&\equiv (1-q)\sum_{k=1}^{\lfloor n/2\rfloor}\frac{(-1)^kq^{k(k+1)/2}}{1-q^{k}} {n-k\brack k}\notag\\[10pt]
&\equiv (q-1)\sum_{k=1}^{\lfloor n/2\rfloor}\frac{(-1)^kq^{k(k-1)/2}}{1-q^{n-k}} {n-k\brack k}\pmod{\Phi_n(q)}.\label{b-3}
\end{align}

We can rewrite \eqref{b-1} as
\begin{align}
\sum_{k=1}^{\lfloor n/2\rfloor}\frac{(-1)^kq^{k(k-1)/2}}{1-q^{n-k}} {n-k\brack k}
=\frac{\mathcal{R}_n(q)-1}{1-q^n}.\label{new-7}
\end{align}
Substituting \eqref{new-7} into the right-hand side of \eqref{b-3} gives
\begin{align}
\sum_{k=1}^{\lfloor n/2\rfloor}\frac{q^{-k(k-1)}}{[2k]_q} {2k\brack k}\equiv
\frac{(1-q)(1-\mathcal{R}_n(q))}{1-q^n}\pmod{\Phi_n(q)},\label{new-8}
\end{align}
as desired.
\qed

\section{Proof of Theorem \ref{t-1}}
For $1\le k \le \lfloor n/2 \rfloor$, we have
\begin{align*}
{n\brack k}&=\frac{(1-q^n)(1-q^{n-1})\dots (1-q^{n-k+1})}{(1-q)\cdots(1-q^k)}\\[10pt]
&\equiv \frac{(1-q^n)(1-q^{-1})\dots (1-q^{-k+1})}{(1-q)\cdots(1-q^k)}\\[10pt]
&=\frac{(-1)^{k-1}q^{-k(k-1)/2}(1-q^n)}{1-q^k}\\[10pt]
&\equiv \frac{(-1)^{k}q^{-k(k+1)/2}(1-q^n)}{1-q^{n-k}}\pmod{\Phi_n(q)^2},
\end{align*}
where we have used the fact that $1-q^n\equiv 0\pmod{\Phi_n(q)}$.
It follows from \eqref{new-2}, \eqref{new-7} and the above that
\begin{align*}
\left({n\choose 0}\right)_q&=\sum_{k=0}^{\lfloor n/2 \rfloor}
q^{k^2}{n \brack k} { n-k \brack k}\notag\\[10pt]
&=1+\sum_{k=1}^{\lfloor n/2 \rfloor}
q^{k^2}{n \brack k} { n-k \brack k}\notag\\[10pt]
&\equiv 1+(1-q^n)\sum_{k=1}^{\lfloor n/2 \rfloor}
 \frac{(-1)^{k}q^{k(k-1)/2}}{1-q^{n-k}}{ n-k \brack k}\pmod{\Phi_n(q)^2}\\[10pt]
 &=\mathcal{R}_n(q).
\end{align*}
This completes the proof of Theorem \ref{t-1}.

\section{Proof of Theorem \ref{t-2}}
Note that for $1\le k\le \lfloor n/2\rfloor$, we have
\begin{align}
{2n \brack k}&=\frac{(1-q^{2n})(1-q^{2n-1})\cdots(1-q^{2n+1-k})}{(1-q)(1-q^2)\cdots(1-q^k)}\notag\\[10pt]
&\equiv \frac{2(1-q^{n})(1-q^{-1})\cdots(1-q^{-k+1})}{(1-q)(1-q^2)\cdots(1-q^k)}\pmod{\Phi_n(q)^2}\notag\\[10pt]
&=\frac{2(-1)^{k-1}(1-q^{n})q^{-k(k-1)/2}}{1-q^k},\label{c-1}
\end{align}
and
\begin{align}
{2n-k \brack n+k}&={2n\brack n}\frac{(1-q^{n-2k+1})\cdots(1-q^n)}{(1-q^{n+1})\cdots (1-q^{n+k})(1-q^{2n-k+1})\cdots(1-q^{2n})}\notag\\[10pt]
&\equiv \frac{(-1)^k q^{-k(3k-1)/2}}{2}{2n\brack n}{2k-1\brack k}\pmod{\Phi_n(q)},\label{c-2}
\end{align}
where we have utilized the fact $q^n\equiv 1\pmod{\Phi_n(q)}$.

It follows from \eqref{new-2}, \eqref{c-1} and \eqref{c-2} that
\begin{align}
\left(2n\choose n\right)_q&=\sum_{k=0}^{\lfloor n/2\rfloor}
q^{k(k+n)}{2n \brack k} { 2n-k \brack n+k}\notag\\[10pt]
&={2n\brack n}+\sum_{k=1}^{\lfloor n/2\rfloor}
q^{k(k+n)}{2n \brack k} { 2n-k \brack n+k}\notag\\[10pt]
&\equiv {2n\brack n}-(1-q^{n}){2n\brack n}\sum_{k=1}^{\lfloor n/2\rfloor}
\frac{q^{-k(k-1)}}{1-q^k}{2k-1\brack k}\pmod{\Phi_n(q)^2}\notag\\[10pt]
&={2n\brack n}-\frac{1-q^{n}}{1-q}{2n\brack n}\sum_{k=1}^{\lfloor n/2\rfloor}
\frac{q^{-k(k-1)}}{[2k]_q}{2k\brack k}.\label{c-3}
\end{align}

From \eqref{a-3}, we deduce that
\begin{align}
{2n\brack n}&\equiv 1+q^{n^2}\notag\\[10pt]
&=2-(1-q^n)(1+q^n+q^{2n}+\cdots+q^{(n-1)n})\notag\\[10pt]
&\equiv 2-n(1-q^n)\pmod{\Phi_n(q)^2}.\label{c-5}
\end{align}

Finally, substituting \eqref{a-1} and \eqref{c-5} into the right-hand side of \eqref{c-3} gives
\begin{align*}
\left(2n\choose n\right)_q&\equiv
2\mathcal{R}_n(q)-n(1-q^n)\pmod{\Phi_n(q)^2},
\end{align*}
as desired.

\vskip 5mm \noindent{\bf Acknowledgments.}
The author would like to thank Ofir Gorodetsky for discussions on $q$-trinomial coefficients and useful suggestions regarding the paper. This work was supported by the National Natural Science Foundation of China (grant 12171370).


\begin{thebibliography}{99}

\small \setlength{\itemsep}{-.8mm}

\bibitem{andrews-jams-1990}G.E. Andrews,
Euler's ``exemplum memorabile inductionis fallacis'' and $q$-trinomial coefficients,
J. Amer. Math. Soc. 3 (1990), 653--669.

\bibitem{andrews-dm-1999}G.E. Andrews, $q$-Analogs of the binomial coefficient congruences of Babbage, Wolstenholme and Glaisher, Discrete Math. 204 (1999), 15--25.

\bibitem{ab-1987}G.E. Andrews and R.J. Baxter, Lattice gas generalization of the
hard hexagon model III: $q$-trinomial coefficients, J. Stat. Phys. 47 (1987), 297--330.

\bibitem{al-integers-2020}M. Apagodu and J.-C. Liu, Congruence properties for the trinomial coefficients, Integers 20 (2020), Art. 38.

\bibitem{babbage-epj-1819}C. Babbage, Demonstration of a theorem relating to prime numbers, Edinburgh Philosophical J. 1 (1819), 46--49.

\bibitem{bachraoui-rj-2021}M. El Bachraoui, On supercongruences for truncated sums of squares of basic hypergeometric series, Ramanujan J. 54 (2021), 415--426.

\bibitem{goro-ijnt-2019}O. Gorodetsky, $q$-Congruences, with applications to supercongruences and the cyclic sieving phenomenon, Int. J. Number Theory 15 (2019), 1919--1968.

\bibitem{guo-aam-2020}V.J.W. Guo, $q$-Supercongruences modulo the fourth power of a cyclotomic polynomial via creative microscoping, Adv. Appl. Math. 120 (2020), Art. 102078.

\bibitem{gs-rm-2020}V.J.W. Guo and M.J. Schlosser, A new family of $q$-supercongruences modulo the fourth power of a cyclotomic polynomial, Results Math. 75 (2020), Art. 155.

\bibitem{gz-am-2019}V.J.W. Guo and W. Zudilin, A $q$-microscope for supercongruences, Adv. Math. 346 (2019), 329--358.

\bibitem{lw-2020}L. Li and S.-D. Wang, Proof of a $q$-supercongruence conjectured by Guo and Schlosser, Rev. R. Acad. Cienc. Exactas F\'is. Nat., Ser. A Mat. 114 (2020), Art. 190.

\bibitem{liu-cmj-2017}J.-C. Liu, Some finite generalizations of Euler's pentagonal number theorem, Czechoslovak Math. J. 142 (2017), 525--531.

\bibitem{liu-crm-2020}J.-C. Liu, On a congruence involving $q$-Catalan numbers, C. R. Math. Acad. Sci. Paris 358 (2020), 211--215.


\bibitem{lp-a-2020}J.-C. Liu and F. Petrov, Congruences on sums of $q$-binomial coefficients,
  Adv. in Appl. Math. 116 (2020), Art. 102003.

\bibitem{np-amc-2018}H.-X. Ni and H. Pan, On the lacunary sum of trinomial coefficients, Appl. Math. Comput. 339 (2018), 286--293.

\bibitem{pan-mm-2013}H. Pan, Factors of some lacunary $q$-binomial sums, Monatsh. Math. 172 (2013), 387--398.

\bibitem{sills-b-2018}A.V. Sills, An invitation to the Rogers-Ramanujan identities, CRC Press, 2018.

\bibitem{straub-pams-2019}A. Straub, Supercongruences for polynomial analogs of the Ap\'ery numbers, Proc. Amer. Math. Soc. 147 (2019), 1023--1036.

\bibitem{sunzw-scm-2011}Z.-W. Sun, Super congruences and Euler numbers, Sci. China Math. 54
    (2011), 2509--2535.

\bibitem{sunzw-scm-2014}Z.-W. Sun, Congruences involving generalized central trinomial coefficients, Sci. China Math. 57 (2014), 1375--1400.

\bibitem{tauraso-aam-2012}R. Tauraso, $q$-Analogs of some congruences involving Catalan numbers, Adv. in Appl. Math. 48 (2012), 603--614.

\bibitem{ws-a-2019}C. Wang and Z.-W. Sun, Congruences involving central trinomial coefficients, preprint (2019), arXiv:1910.06850.

\bibitem{WY2}X. Wang and M. Yue, A $q$-analogue of a Dwork-type supercongruence, Bull. Aust. Math. Soc. 103 (2021), 303--310.

\bibitem{wei-jcta-2021}C. Wei, Some $q$-supercongruences modulo the fourth power of a cyclotomic polynomial, J. Combin. Theory Ser. A 182 (2021), Art. 105469.

\bibitem{wolstenholme-qjpam-1862}J. Wolstenholme, On certain properties of prime numbers, Quart. J. Pure Appl. Math. 5 (1862), 35--39.

\bibitem{zudilin-ac-2019}W. Zudilin, Congruences for $q$-binomial coefficients, Ann. Comb. 23 (2019), 1123--1135.

\end{thebibliography}
\end{document}